# A New Approach to Reducing Vector Delay Nonlinear Systems: Boundedness and Stability Analysis

Mark A. Pinsky

**Abstract.** This paper introduces a new method for assessing the boundedness and stability of certain vector nonlinear systems with delays and variable coefficients. The approach is based on developing scalar counterparts to the given vector systems. We prove that the solutions to these scalar nonlinear equations, which also include delays and variable coefficients, provide upper bounds for the norms of solutions to the original vector equations if the history functions for both equations are properly matched. This enables the evaluation of the boundedness and stability characteristics of a vector system by analyzing the abridged dynamics of its scalar counterparts. This assessment can be carried out through straightforward simulations or by applying simplified analytical methods. As a result, we introduce new criteria for boundedness and stability and estimate the radii of the balls that contain history functions stemming bounded or stable solutions for the original vector systems. Finally, we validate our inferences through representative simulations that also assess the accuracy of our approach.

**Keywords.** Boundedness/ Stability of nonlinear delay systems, Boundedness/stability regions, Delay Nonlinear Systems, Comparison principle, Estimation of solutions' norms, Model reduction, Variable Coefficients.

## 1. Introduction

Assessing the boundedness and stability of vector nonlinear time-delay systems (VNDS) with variable coefficients is a challenging problem that has significant implications in various fields such as biology, medicine, economics, physics, and control science. Accordingly, this problem has been extensively studied in numerous monographs [8],[9],[12]-[14],[16],[18],[19],[21],23],[25],[28] and research papers, see also reviews [2],[3],[15],[20] that develop stability criteria for certain types of time-delay systems. Key advancements in this area have been made through extensions of the Lyapunov methodology on time-delay systems, initially developed by N. Krasovskii [24] and Razumikhin [33]. This methodology, along with the application of Halanay's inequality [17], see also [12], [13], and the Bohl–Perron theorem [4],[5],[14], and certain other techniques has led the development of multiple stability criteria for delay systems.

For linear delay systems, the application of the Lyapunov-like techniques have resulted in a computationally efficient method based on linear matrix inequalities, as reviewed in numerous studies, see, e.g., [13],[5],[16], [22],[23] and [34] and additional citation therein. The stability of time-varying delay systems, which sometimes include time-varying delays, has also been undertaken in [5],[8], [9],[10],[26] and [37].

However, the application of the Lyapunov methodology to time-varying nonlinear systems in finite dimensions (ODEs) already has become fairly arduous and often results in overly conservative stability criteria. Extending this approach to infinite dimensional VNDS presents additional challenges, particularly for the systems coupling multiple nonlinear delay equations with variable coefficients. As a result, most current techniques provide conservative stability criteria that rarely define the regions of stability or attraction with adequate accuracy for non-autonomous VNDS.

These limitations are even more pronounced when the Lyapunov-like methodology is used to assess the boundedness of solutions to nonhomogeneous VNDS with variable coefficients, a topic that has been rarely addressed, see, for instance [1]. Although input-to-state stability analysis has addressed some aspects of the behavior of such systems, as highlighted in recent review paper [3], many challenges still remain.

The dynamic behavior of some scalar nonlinear systems with delay and time-varying coefficients has also been studied, providing additional insights into dynamics of these systems, see, e.g., [6] and [7], and additional references therein.

Lastly, we note that the study of finite-time stability, which examines the quantitative and qualitative characteristics of the underlying problems, has been explored for delay systems as well. For further details, see a recent review paper [35] and the additional references cited therein.

This paper presents a novel approach for assessing the behavior of the norms of solutions to nonautonomous VNDS by developing scalar auxiliary counterparts for these systems. Building on our previous studies of similar systems in finite dimensions [29]-[32], we show that solutions to these scalar delay equations provide an upper

_________________________________________________________________________________
Mark A. Pinsky, Department of Mathematics and Statistics, University of Nevada. Reno, Reno NV 89557, USA, email: pinsky@unr.edu.

bound on the evolution of the norms of the original VNDS if the history functions for both systems are appropriately matched.

Consequently, this approach enables the evaluation of the boundedness and stability properties of nonautonomous VNDS by analyzing the dynamics of their scalar counterparts that can be readily simulated or analyzed through simplified reasoning. As a result, we derive new boundedness and stability criteria and estimate the radii of the balls that contain history functions stemming bounded or stable solutions to the original systems. Finally, we validate this approach through representative simulations, which also assess the accuracy of these techniques.

This paper is organized as follows. The next section outlines our notation, some preliminary statements, and defines the underlined system. Section 3 ascribes the reduction technique and its applications to the evaluation of the boundedness/stability of some VNDS with variable coefficients. Section 4 includes a closed-form robust stability criterion. Section 5 discusses the applications of linearized auxiliary equations. Section 6 highlights the results of our simulations, and Section 7 concludes this study and outlines some directions for subsequent research.

## 2. Notation, Mathematical Preliminaries, and Governing Equations

**2.1. Notation**. Firstly, let us recall that symbols $\mathbb{R}$, $\mathbb{R}_{\geq 0}$, $\mathbb{R}_+$ and $\mathbb{R}^n$ stand for the sets of real, non-negative and positive real numbers, and real $n$-dimensional vectors, $\mathbb{N}$ is a set of real positive integers, $\mathbb{R}^{n \times n}$ is a set of $n \times n$-matrices, and $I \in \mathbb{R}^{n \times n}$ is the identity matrix. Next, $C([a,b]; \mathbb{R}^n)$, $C([a,b]; \mathbb{R}_+)$ and $C([a,b]; \mathbb{R}_{\geq 0})$ are the spaces of continuous real functions $\zeta : [a,b] \to \mathbb{R}^n$, $\zeta : [a,b] \to \mathbb{R}_+$ or $\zeta : [a,b] \to \mathbb{R}_{\geq 0}$, respectively, with the supremum norm $\|\zeta\| := \sup_{t \in [a,b]} |\zeta(t)|$, where $|\cdot|$ stands for the Euclidean norm of a vector or the induced norm of a matrix, $a < b \in \mathbb{R}$ and $b$ can be infinity. In turn, $\dot{x}(t) := D^+ x(t)$, where $D^+$ is the upper righthand derivative in $t$. Also note that $|x|_\infty = \sup_{i=1,\ldots,n}(|x_i|)$ and $|x|_1 = \sum_{i=1}^{n} |x_i|$, $\forall x \in \mathbb{R}^n$.

**2.2. Preliminaries.** The solutions to a scalar differential equation are subject to the comparison statements that are extended in this section on the solutions to some scalar delay equations using method of steps and mathematical induction. For this reason, we acknowledge two pertinent comparison statements for scalar ODEs [36].

**Lemma 1**. Consider a scalar equation,

$$\dot{u} = f(t,u), \ u(t_0) = u_0 \tag{2.1}$$

where function $f \in C([t_0, \infty) \times \mathbb{R}; \mathbb{R})$ is locally Lipschitz in $u$, $\forall u \in J \in \mathbb{R}$. Next, we presume that (2.1) admits a unique solution $u(t, u_0) \in J$, $\forall t \in [t_0, \infty)$ and write a coupling scalar inequality,

$$\dot{v} \leq f(t,v), \ v(t_0) = v_0 \leq u_0 \tag{2.2}$$

where $v(t, v_0) \in J$, $\forall t \geq t_0$ is a solution of (2.2). Then

$$v(t, v_0) \leq u(t, u_0), \ \forall t \in [t_0, t_*] \tag{2.3}$$

where $t_*$ can be infinity. The next lemma is often derived through the application of Lemma 1 [36].

**Lemma 2**. Consider two scalar equations,

$$\dot{u}_1 = f_1(t, u_1), \ u_1(t_0) = u_{10},$$
$$\dot{u}_2 = f_2(t, u_2), \ u_2(t_0) = u_{20},$$

where functions $f_i \in C([t_0, t_*] \times \mathbb{R}; \mathbb{R})$, $i = 1, 2$ are locally Lipschitz in $u_i$, $\forall u_i \in J_* \subset \mathbb{R}$, and $t_*$ can be infinity. Presume that both above equations admit unique solutions, $u_i(t, u_{i0}) \in J_*, \forall t \in [t_0, t_*]$, $i = 1, 2$ and that $f_1(t, u) \leq f_2(t, u)$, $\forall u \in J_*$, $\forall t \in [t_0, t_*]$ and $u_{10} \leq u_{20}$. Then $u_1(t, u_{10}) \leq u_2(t, u_{20})$, $\forall t \in [t_0, t_*]$.

Subsequently, using the method of steps, we extend these statements to a scalar functional differential equation and the corresponding inequality with multiple variable delays $h_i(t)$ that are subject to the following conditions:

$$h_i \in C([t_0,\infty);\mathbb{R}_+), \ \max_i \sup_{\forall t \geq t_0} h_i(t) = \bar{h} < \infty, \ \min_i \inf_{\forall t \geq t_0} h_i(t) \geq \underline{h} > 0, \ i=1,...,m \quad (2.4)$$

that ensure the application of time stepping with, e.g., step equals $\underline{h}$. In turn, we assume that the endpoints of the consecutive time-steps are set at $s_0 = t_0$, $s_k = t_0 + k\bar{h}$, $k \in \mathbb{N}$. This leads to the following statement.

**Lemma 3**. Denote a scalar delay differential equation and the coupling inequality as follows:

$$\dot{u} = f(t, u(t), u(t-h_1(t)),...,u(t-h_m(t))), \ \forall t \geq t_0, \ u(t) = \varphi(t), \ \forall t \in [t_0 - \bar{h}, t_0] \quad (2.5)$$

$$\dot{v} \leq f(t, v(t), v(t-h_1(t)),...,v(t-h_m(t))), \ \forall t \geq t_0, \ v(t) = \phi(t), \ \forall t \in [t_0 - \bar{h}, t_0] \quad (2.6)$$

where a continuous function $f \in \mathbb{R}$ is locally Lipschitz in the second variable $\forall u \in J_\times \subset \mathbb{R}$ and non-decreasing in all other variables starting from the third, $\phi, \varphi \in C([t_0 - \bar{h}, t_0]; \mathbb{R})$, and $h_i(t)$ are subject to (2.4). Additionally, assume that equation (2.5) admits a unique solution, $u(t,\varphi) \in J_\times$, $\forall t \geq t_0$, $\forall \|\varphi\| \leq \bar{\varphi} \in \mathbb{R}_+$. Then

$$v(t, \phi) \leq u(t, \varphi), \ \forall t \in [t_0, \infty) \quad (2.7)$$

where $v(t, \phi)$ is a solution to (2.6).

**Proof**. Let us set that $u_k(t) := u(t, \varphi)$ and $v_k(t) := v(t, \phi)$, $\forall t \in [s_{k-1}, s_k]$, $k \geq 1$. Then on the first-time step (2.5) and (2.6) yield that

$$\dot{u}_1 = f(t, u_1(t), \varphi(t-h_1(t)),...,\varphi(t-h_m(t))) = f_1^1(t, u_1), \ \forall t \in [t_0, t_0+s_1]$$
$$\dot{v}_1(t) \leq f(t, v_1(t), \phi(t-h_1(t)),...,\phi(t-h_m(t))) = f_2^1(t, v_1), \ \forall t \in [t_0, t_0+s_1] \quad (2.8)$$

Since $\varphi(t) \geq \phi(t)$, $\forall t \in [t_0, t_0+s_1]$ and $f$ is nondecreasing in the variables that include delay, we infer that $f_1^1(t,u) \geq f_2^1(t,u)$, $\forall t \in [t_0, t_0+s_1]$, $\forall u \in J_\times$. In turn, using the last inequality, we can write the second equation in (2.8) as $\dot{v}_1(t) \leq f_1^1(t, v_1)$, $\forall t \in [t_0, t_0+s_1]$. Then, the application of Lemma 1 to (2.8) grants that $v_1(t) \leq u_1(t)$, $\forall t \in [t_0, s_1]$. Assume, in turn, that $v_k(t) \leq u_k(t)$, $\forall t \in [s_{k-1}, s_k]$, $k > 1$. Then on the next time step, (2.5) and (2.6) can be written as follows,

$$\dot{u}_{k+1} = f(t, u_{k+1}(t), u_k(t-h_1(t)),...,u_k(t-h_m(t))) = f_1^{k+1}(t, u_{k+1}), \ \forall t \in [s_k, s_{k+1}]$$
$$\dot{v}_{k+1}(t) \leq f(t, v_{k+1}(t), v_k(t-h_1(t)),...,v_k(t-h_m(t))) = f_2^{k+1}(t, v_{k+1}), \ \forall t \in [s_k, s_{k+1}] \quad (2.9)$$

where $f_1^{k+1}(t,u) \geq f_2^{k+1}(t,u)$, $\forall t \in [s_k, s_{k+1}]$, $\forall u \in J_\times$ due to both the assumption on function $f$ and the induction hypothesis. Then the application of the latter relation to (2.9) implies that,

$$\dot{v}_{k+1}(t) \leq f(t, v_{k+1}, u_k(t-h_1(t)),...,u_k(t-h_m(t))) = f_1^{k+1}(t, v_{k+1}), \ \forall t \in [s_k, s_{k+1}] \quad (2.10)$$

Thus, the first equation (2.9) and the equation (2.10) are matched with equations (2.1) and (2.2), respectively, prompting the application of Lemma 1 □

In turn, the extension of Lemma 2 to scalar delay equations is enabled by the following statement.

**Lemma 4**. Let us write two scalar equations with delay in the following form,

$$\dot{u}_1 = f_1(t, u_1(t), u_1(t-h_1(t)),...,u_1(t-h_m(t))), \ \forall t \in [t_0, t_*], \ u_1(t) = \varphi(t), \ \forall t \in [t_0 - \bar{h}, t_0]$$
$$\dot{u}_2 = f_2(t, u_2(t), u_2(t-h_1(t)),...,u_2(t-h_m(t))), \ \forall t \in [t_0, t_*], \ u_2(t) = \phi(t), \ \forall t \in [t_0 - \bar{h}, t_0] \quad (2.11)$$

where continuous functions $f_i \in \mathbb{R}$, $i=1,2$ are locally Lipschitz in the second variables, $\forall u_i \in J_+ \subset \mathbb{R}$, $i=1,2$, $\phi, \varphi \in C([t_0 - \bar{h}, t_0]; \mathbb{R})$, and $t_*$ can be infinity. Next, we assume that both the above equations admit the unique solutions $u_i(t, \varphi) \in J_+$, $\forall t \in [t_0, t_*]$, $\forall \|\varphi\|, \forall \|\phi\| \leq \bar{\varphi} \in \mathbb{R}_+$ and that,

$$f_1(t, x_1,...,x_{m+1}) \leq f_2(t, x_1,...,x_{m+1}), \ \forall t \in [t_0, t_*], \ \forall x_i \in \mathbb{R}, \ i=2,...,m+1,$$
$$\varphi(t) \leq \phi(t), \ \forall t \in [t_0 - \bar{h}, t_0] \quad (2.12)$$

Then $u_1(t,\varphi) \le u_2(t,\phi)$, $\forall t \in [t_0, t_*]$.

**Proof**. The proof of this statement is similar to the proof of the previous lemma. In fact, let us set that $u_1^k(t) := u_1(t,\varphi)$, $u_2^k(t) := u_2(t,\phi)$, $\forall t \in [s_{k-1}, s_k]$, $k \ge 1$. Then we initially write that

$$\dot{u}_1^1 = f_1(t, u_1^1(t), \varphi(t-h_1(t)),...,\varphi(t-h_m(t))) = f_1^1(t, u_1^1), \forall t \in [t_0, t_0+s_1],$$
$$\dot{u}_2^1 = f_2(t, u_2^1(t), \phi(t-h_1(t)),...,\phi(t-h_m(t))) = f_2^1(t, u_2^1), \forall t \in [t_0, t_0+s_1]$$

Consequently, (2.12) implies that $f_1^1(t,u) \le f_2^1(t,u)$, $\forall u \in J_+$, $\forall t \in [t_0, t_0+s_1]$. Then Lemma 2 yields that $u_1^1(t) \le u_2^1(t)$, $\forall t \in [t_0, t_0+s_1]$. In turn, we assume that $u_1^k(t) \le u_2^k(t)$, $\forall t \in [s_{k-1}, s_k]$ and show that the latter inequality can be extended for $k := k+1$. For this sake we write that

$$\dot{u}_1^{k+1} = f_1(t, u_1^{k+1}(t), u_1^k(t-h_1(t)),...,u_1^k(t-h_m(t))) = f_1^{k+1}(t, u_1^{k+1}), \forall t \in [s_k, s_{k+1}],$$
$$\dot{u}_2^{k+1} = f_2(t, u_2^{k+1}(t), u_2^k(t-h_1(t)),...,u_2^k(t-h_m(t))) = f_2^{k+1}(t, u_2^{k+1}), \forall t \in [s_k, s_{k+1}] \quad (2.13)$$

and conclude, due to both (2.12) and the induction's conjecture that, $f_1^{k+1}(t,u) \le f_2^{k+1}(t,u)$, $\forall u \in J_+$, $\forall t \in [s_k, s_{k+1}]$ which together with Lemma 2 ensures this statement □

2.3 **Underlined Equation and Definitions**. In the sequel, we are going to study the behavior of solutions of the following vector nonlinear equation with multiple variable delays and variable coefficients,

$$\dot{x} = A(t)x + f(t, x(t), x(t-h_1(t)),...,x(t-h_m(t))) + F(t), \forall t \ge t_0, x(t,\varphi) = \varphi(t), \forall t \in [t_0-\bar{h}, t_0] \quad (2.14)$$

where $x \in N \subset \mathbb{R}^n$, $0 \in N$, a function $f \in \mathbb{R}^n$ is continuous in all variables and locally Lipschitz in the second one, $f(t,0) = 0$, $\varphi \in C([t_0-\bar{h}, t_0]; \mathbb{R}^n)$, $\|\varphi\| := \sup_{t \in [t_0-\bar{h}, t_0]} |\varphi(t)|$, continuous functions $h_i(t)$ are defined by (2.4), matrix $A \in C([t_0, \infty); \mathbb{R}^{n \times n})$, $F(t) = F_0 e(t)$, $e \in C([t_0, \infty); \mathbb{R}^n)$, $\|e\| = \sup_{t \ge t_0} |e(t)| = 1$, and $F_0 \in \mathbb{R}_{\ge 0}$.

Furthermore, we presume that the initial problem (2.14) assumes a unique solution $\forall t \ge t_0$ and $\forall \|\varphi\| \le \bar{\phi} > 0$.

Under these assumptions, the right side of (2.14) is a continuous vector function in all variables, which implies that a solution to this equation, $x(t, t_0, \varphi)$ is continuous and continuously differentiable in $t$, $\forall t \ge t_0$ and, thus, is bounded on any finite time interval. Consequently, this paper primarily focuses on the behavior of solutions to (2.14) as $t \to \infty$.

Note that in the sequel we shall employ the abridged notation for the solutions to (2.14) as follows $x(t,\varphi) := x(t, t_0, \varphi)$, $\forall t \ge t_0$.

To simplify further referencing, we also acknowledge the homogeneous counterpart of (2.14),

$$\dot{x} = A(t)x + f(t, x(t), x(t-h_1(t)),...,x(t-h_m(t))), \forall t \ge t_0, x(t,\varphi) = \varphi(t), \forall t \in [t_0-\bar{h}, t_0] \quad (2.15)$$

and the following linear equation,

$$\dot{x} = A(t)x, \forall t \ge t_0, x(t_0, \varphi) = \varphi(t_0) \in \mathbb{R}^n \quad (2.16)$$

Subsequently, we assume that (2.15) possesses a unique solution $\forall \|\varphi\| \le \bar{\phi}$, $\forall t \ge t_0$ and write the solution to (2.16) as follows, $x(t) = W(t, t_0) \varphi(t_0)$, where $W(t, t_0) = w(t) w^{-1}(t_0)$ is the transition (Cauchy) matrix and $w(t)$ is a fundamental solution matrix for (2.16). In the sequel, we shall assume that $|w(t_0)| = 1$.

Next, we present the standard definitions of the boundedness/stability of solutions to either equations (2.14) or (2.15) which will be used in the sequel of this paper, see, e.g.,[23].

**Definition 1**. Assume that $\varphi \in C([t_0-\bar{h}, t_0]; \mathbb{R}^n)$ and (2.15) admits a unique solution $\forall \|\varphi\| \le \bar{\phi} > 0$. Then the trivial solution of equation (2.15) is called:

1) stable for the set value of $t_0$ if $\forall \varepsilon \in \mathbb{R}_+$, $\exists \delta_1(t_0, \varepsilon) \in \mathbb{R}_+$ such that $\forall \|\varphi\| < \delta_1(t_0, \varepsilon), |x(t, t_0, \varphi)| < \varepsilon, \forall t \ge t_0$. Otherwise, the trivial solution is called unstable.

2) uniformly stable if in the above definition $\delta_1(t_0,\varepsilon) = \delta_2(\varepsilon)$.

3) asymptotically stable if it is stable for given value of $t_0$ and $\exists \delta_3(t_0) \in \mathbb{R}_+$ such that $\lim_{t\to\infty}|x(t,t_0,\varphi)| = 0$, $\forall \|\varphi\| < \delta_3(t_0)$.

4) uniformly asymptotically stable if it is uniformly stable and in the previous definition $\delta_3(t_0) = \delta_4 = const$.

5) uniformly exponentially stable if $\exists \delta_5 \in \mathbb{R}_+$ and $\exists c_i \in \mathbb{R}_+$, $i=1,2$ such that,
$$|x(t,t_0,\varphi)| \leq c_1 \|\varphi\| \exp(-c_2(t-t_0)), \forall \|\varphi\| \leq \delta_5, \forall t \geq t_0$$

**Definition 2.** Let (2.14) admits a unique solution $\forall \|\varphi\| \leq \bar{\phi}$, $\varphi \in C([t_0-\bar{h},t_0];\mathbb{R}^n)$. A solution to equation (2.14) is called:

6) bounded for the set value $t_0$ if $\exists \delta_6(t_0) \in \mathbb{R}_+$, $\exists F_*(t_0) \in \mathbb{R}_{\geq 0}$ and $\exists \varepsilon_*(\delta_6,F_*) \in \mathbb{R}_+$ such that $|x(t,t_0,\varphi)| < \varepsilon_*$ $\forall t \geq t_0$, $\forall \|\varphi\| < \delta_6(t_0)$ and $\forall F_0 \leq F_*(t_0)$.

7) uniformly bounded if both $\delta_6(t_0) = \delta_7$ and $F_*(t_0) = \hat{F}$.

Furthermore, let $R_i \equiv \sup \delta_i$, $i=1,...,7$ be the superior value of $\delta_i$ for which either the $i$-th condition in Definition 1, or conditions 6 or 7 in Definition 2 hold. This implies that $R_1 = R_1(t_0)$, $R_3 = R_3(t_0)$, $R_6 = R_6(t_0,F_0)$, $\forall F_0 \leq F_*(t_0)$ and $R_7 = R_7(F_0)$, $\forall F_0 \leq \hat{F}$. Next, let $B_{R_i} := \{\varphi \in C([t_0-s_1,t_0],\mathbb{R}^n) : \|\varphi\| \leq R_i\}$, $i=1,...,7$ be a ball with the radius $R_i$ that is centered at zero. Clearly, these balls encompass history functions that originate bounded/ stable solutions to (2.14)/(2.15), respectively.

Next, we recall a definition of robust stability that was termed stability under persistent perturbations by N. Krasovskii [24], pp.161-164]. Consider a perturb equation to (2.15) in the following form,
$$\dot{x} = A(t)x + f(t,x,x(t-h_1(t)),...,x(t-h_m(t))) + R(t,x,x(t-h_1^*(t)),...,x(t-h_1^*(t))), \forall t \geq t_0,$$
$$x(t) = \varphi(t), \forall t \in [t_0-\bar{h},t_0] \quad (2.17)$$

where, in addition to the conditions set above for the components of (2.17) that are included in (2.15), we also assume that (I) function $f(t,\chi_1,...,\chi_{m+1})$ is Lipschitz continuous in $\chi_i \in \mathbb{R}^n$, $i=1,...,m+1$, (II) $f(t,0)=0$, (III) $R \in \mathbb{R}^n$ is a continuous function in all its variables but $R(t,0)$ might not be zero, and (IV) $h_i^* \in C([t_0,\infty);\mathbb{R}_+)$ comply with the conditions (2.4).

**Definition 3.** The trivial solution to (2.15) is called robustly stable if $\forall \varepsilon \in \mathbb{R}_+$, $\exists \Delta_i(\varepsilon) \in \mathbb{R}_+$, $i=1,2,3$ such that $|x(t,t_0,\varphi)|_\infty \leq \varepsilon$, $\forall t \geq t_0$ if $|R(t,\chi_1,...,\chi_{m+1})|_\infty < \Delta_1(\varepsilon)$, $\forall t \geq t_0$, $\forall |\chi_i|_\infty < \varepsilon$, $i=1,...,m+1$, $\|\varphi\| < \Delta_2(\varepsilon)$, $|h_i(t)-h_i^*(t)| < \Delta_3(\varepsilon)$, $\forall t \geq t_0$, $i=1,...,m$, where $x(t,t_0,\varphi)$ is a solution to (2.17).

Clearly, this definition commences the robust stability of (2.15) through the stability of defined above perturbed equation (2.17).

**Statement** (N. Krasovskii [24], p.162). Assume that conditions (I)-(IV) hold and that the trivial solution of (2.15) is uniformly asymptotically stable. Then it is robustly stable as well.

Lastly, we extend a basic definition of finite-time stability (FTS) that was used for ODEs in [11] onto the systems with delay.

**Definition 4.** Equation (2.14) or (2.15) is called
  a) finite-time stable (FTS) with respect to positive numbers $\alpha,\beta,T$, $\alpha < \beta$ if the condition $\|\varphi\| < \alpha$ implies that $|x(t,t_0,\varphi)| < \beta(\alpha,T)$, $\forall t \in [t_0,t_0+T]$.
  b) finite-time contractively stable (FTCS) with respect to positive numbers $\alpha,\beta,\gamma,T$, $\alpha < \gamma$ if it is FTS with respect of $\alpha,\beta,T$ and $\exists t_1 \in (t_0,t_0+T)$ such that $|x(t,t_0,\varphi)| < \gamma(\alpha,T)$, $\forall t \in [t_1,t_0+T]$.

**3. Scalar Delay Auxiliary Equations**

This section derives a scalar delay auxiliary differential equation which prompts the estimation from the above of the time histories of the norms of solutions to some vector nonlinear systems with variable coefficients and delays.

First, using the variation of parameters, we write (2.14) as follows,

$$x(t,\varphi) = w(t)w^{-1}(t_0)\varphi(t_0) + w(t)\int_{t_0}^{t} w^{-1}(\tau)\Big(f\big(t, x(\tau), x(\tau-h_1(\tau)), ..., x(\tau-h_m(\tau))\big) + F(\tau)\Big)d\tau, \quad \forall t \geq t_0,$$

$$x(t,\varphi) = \varphi(t), \quad \forall t \in \left[t_0 - \bar{h}, t_0\right]$$

The last equation implies that,

$$|x(t,\varphi)| = \left| w(t)w^{-1}(t_0)\varphi(t_0) + w(t)\int_{t_0}^{t} w^{-1}(\tau)\Big(f\big(t, x(\tau), x(\tau-h_1(\tau)), ..., x(\tau-h_m(\tau))\big) + F(\tau)\Big)d\tau \right|, \quad \forall t \geq t_0, \quad (3.1)$$

$$|x(t,\varphi)| = |\varphi(t)|, \quad \forall t \in \left[t_0 - \bar{h}, t_0\right]$$

Then the application of the standard norm's inequalities embraces the following equation,

$$|X_1(t,\varphi)| = |w(t)|\left|w^{-1}(t_0)\varphi(t_0)\right| +$$

$$|w(t)|\int_{t_0}^{t}\left|w^{-1}(\tau)\right|\left|f\big(\tau, X_1(\tau), X_1(\tau-h_1(\tau)), ..., X_1(\tau-h_m(\tau))\big) + F(\tau)\right|d\tau, \quad \forall t \geq t_0, \quad (3.2)$$

$$|X_1(t,\varphi)| = |\varphi(t)|, \quad \forall t \in \left[t_0 - \bar{h}, t_0\right]$$

Comparing (3.1) and (3.2) yields that $|x(t,\varphi)| \leq |X_1(t,\varphi)|$, $\forall t \geq t_0$, where $X_1(t,\varphi) \in \mathbb{R}^n$.

To write (3.2) in a more tractable form, we introduce a nonlinear extension of the Lipschitz continuity condition as follows,

$$|f(t,\chi_1,...,\chi_{m+1})| \leq L(t,|\chi_1|,...,|\chi_{m+1}|), \quad \forall \chi = [\chi_1,...,\chi_{m+1}]^T \in \Omega \in \mathbb{R}^{n(m+1)}, \; \chi_i \in \mathbb{R}^n, \; i = 2,...,m+1, \quad (3.3)$$

where $\Omega$ is a compact subset of $\mathbb{R}^{n(m+1)}$ containing zero, $f \in C\big([t_0,\infty) \times \mathbb{R}^{n(m+1)}; \mathbb{R}^n\big)$ and $L \in C\big([t_0,\infty) \times \mathbb{R}^{m+1}_{\geq 0}; \mathbb{R}_{\geq 0}\big)$ is a scalar continuous function and $L(t,0) = 0$.

We illustrate in Appendix how to define such scalar functions in a close form if $f(t,\chi_1,...,\chi_{m+1})$ is a polynomial or power series in $\chi_1,...,\chi_{m+1}$ with, for instance, a bounded in $\Omega$ error term. In the former case $\Omega \equiv \mathbb{R}^{n(m+1)}$ and this condition is also holds if the error term of power series is bounded in $\mathbb{R}^{n(m+1)}$. Also, $L$ turns out to be a linear function in $|\chi_i|$ if $f$ is a linear function in the coresponding variables.

Note that in the remainder of this paper we assume for simplicity that $\Omega \equiv \mathbb{R}^{n(m+1)}$.

Next, the utility of (3.3) in (3.2) brings the following equation,

$$X_2(t,\tilde{\varphi}) = |w(t)|\left|w^{-1}(t_0)\varphi(t_0)\right| +$$

$$|w(t)|\int_{t_0}^{t}\left|w^{-1}(\tau)\right|\Big(L\big(\tau, X_2(\tau), X_2(\tau-h_1(\tau)), ..., X_2(\tau-h_m(\tau))\big) + |F(\tau)|\Big)d\tau, \quad \forall t \geq t_0 \quad (3.4)$$

$$X_2(t,\tilde{\varphi}) = |\varphi(t)|, \quad \forall t \in \left[t_0 - \bar{h}, t_0\right]$$

where $X_2(t) \in \mathbb{R}_{\geq 0}$. Hence, (3.4) implies that $|x(t,\varphi)| \leq X_2(t,\bar{\varphi})$, $\forall t \geq t_0$.

Next, we match (3.4) with the solution to the initial problem of the following scalar delay differential equation,

$$\dot{X}_3 = p(t)X_3 + c(t)\Big(L\big(t, X_3(t), X_3(t-h_1(t)), ..., X_3(t-h_m(t))\big) + |F(t)|\Big), \quad \forall t \geq t_0$$

$$X_3(t,\phi) = \phi(t), \quad \forall t \in \left[t_0 - \bar{h}, t_0\right] \quad (3.5)$$

where, $X_3(t,\phi) \in \mathbb{R}_{\geq 0}$, $p:[t_0,\infty) \to \mathbb{R}$, $c:[t_0,\infty) \to [1,\infty]$ and $\phi \in C\big(\left[t_0 - \bar{h}, t_0\right]; \mathbb{R}_{\geq 0}\big)$.

To define functions $p(t)$, $c(t)$ and $\phi(t)$, we turn (3.5) into its integral form using the variation of parameters,

$$X_3(t) = e^{d(t)}\left(\phi(t_0) + \int_{t_0}^{t} e^{-d(\tau)} c(\tau) q\left(\tau, X_3(\tau), X_3(\tau - h_1(\tau)), \ldots, X_3(\tau - h_m(\tau))\right) d\tau\right), \forall t \geq t_0 \quad (3.6)$$

$$X_{\backslash 3}(t, \phi) = \phi(t), \forall t \in [t_0 - \bar{h}, t_0]$$

where, $d(t) = \int_{t_0}^{t} p(s) ds$ and $q = L(t, X_3(t), X_3(t - h_1(t)), \ldots, X_3(t - h_m(t))) + |F(t)|$.

To determine $p(t)$ and $c(t)$, we shall match the right sides of (3.6) and (3.4). Matching the first additions in these expressions, i.e., $|w(t)||w^{-1}(t_0)\varphi(t_0)|$ and $e^{d(t)}\phi(t_0)$, returns that,

$$|w(t)| = \exp\left(\int_{t_0}^{t} p(s) ds\right) \quad (3.7)$$

$$\phi(t_0) = |w^{-1}(t_0)\varphi(t_0)| \quad (3.8)$$

In turn, to match the second additions in the right side of (3.6) and (3.4), we firstly multiply and divide the latter function by $|w(t)|$ and use that, due to (3.7), $e^{-d(t)} = 1/|w(t)|$, which returns,

$$\int_{t_0}^{t} e^{-d(\tau)} c(\tau) q(\tau, x, x(\tau - h_1(\tau)), \ldots, x(t - h_m(\tau))) d\tau =$$

$$\int_{t_0}^{t} (|w(\tau)||w^{-1}(\tau)|) q(\tau, x, x(\tau - h_1(\tau)), \ldots, x(t - h_m(\tau))) / |w(\tau)| d\tau$$

The last relation yields that

$$c(t) = |w(t)||w^{-1}(t)|$$

is the running condition number of $w(t)$. Let us recall that $0 < |w^{-1}(t)| < \infty$, $\forall t \geq t_0$ which implies that $c(t) < \infty$, $\forall t \geq t_0$. Additionally, continuity of $A(t)$ implies continuity of both $p(t)$ and $c(t)$.

Lastly, from (3.7) we get that,

$$p(t) = d(\ln|w(t)|)/dt \quad (3.9)$$

Next, matching the initial functions in (3.4) and (3.6) yields that $\phi(t) = |\varphi(t)|$, $\forall t \in [t_0 - \bar{h}, t_0]$ which together with (3.8) brings the following condition, $|w^{-1}(t_0)\varphi(t_0)| = |\varphi(t_0)|$ which implies that $w(t_0)$ is the orthonormal matrix. For instance, we can set $w(t_0) = I$ to meet our previous condition that $|w(t_0)| = 1$.

Finally, we write (3.5) as follows,

$$\dot{y} = p(t) y + c(t)\left(L(t, y(t), y(t - h_1(t)), \ldots, y(t - h_m(t))) + |F(t)|\right) \quad (3.10)$$

$$y(t, \phi) = |\varphi(t)|, \forall t \in [t_0 - \bar{h}, t_0]$$

where $y(t) \in \mathbb{R}_{\geq 0}$. This embraces the following statement.

**Theorem 1**. Assume that $f(t, \chi_1, \ldots, \chi_{m+1}) \in \mathbb{R}^n$, $\chi_i \in \mathbb{R}^n$, $i = 2, \ldots, m+1$ is a continuous function in all variables and locally Lipschitz in the second one, $f(t, 0) = 0$, $\varphi \in C([t_0 - \bar{h}, t_0]; \mathbb{R}^n)$, scalar and continuous functions $h_i(t)$ are defined by (2.4), $A \in C([t_0, \infty); \mathbb{R}^{n \times n})$, $F(t) = F_0 e(t)$, $e \in C([t_0, \infty); \mathbb{R}^n)$, $\|e\| = 1$, $F_0 \in \mathbb{R}_{\geq 0}$, a scalar function $L(t, \zeta_1, \ldots, \zeta_{m+1}) \in \mathbb{R}_{\geq 0}$, $\zeta_i \in \mathbb{R}_{\geq 0}$, $i = 2, \ldots, m+1$ is continuous in all variables and locally Lipschitz in the second variable, $L(t, 0) = 0$, inequality (3.3) holds with $\Omega \equiv \mathbb{R}^{n(m+1)}$ and $w(t_0) = I$. Subsequently, we assume that equations (2.14) and (3.10) assume the unique solutions $\forall \|\varphi\| \leq \bar{\varphi} \in \mathbb{R}_+$, $\forall t \geq t_0$. Then

$$|x(t,\varphi)| \leq y(t,\phi), \ \forall t \geq t_0$$
$$\phi(t) = |\varphi(t)|, \ \forall t \in \left[t_0 - \bar{h}, t_0\right] \quad (3.11)$$

where $x(t,\varphi)$ and $y(t,\phi) \geq 0$ are the solutions to (2.14) and (3.10).

**Proof**. In fact, it was shown previously that $|x(t,\varphi)| \leq X_2(t,\tilde{\varphi}), \ \forall t \geq t_0$. Made above assignments of the functions $p(t)$, $c(t)$ and $\phi(t)$ embraces that $X_2(t,\tilde{\varphi}) \equiv X_3(t,\phi) \equiv y(t,\phi), \ \forall t \geq t_0$ if $\tilde{\varphi}(t) = \phi(t) = |\varphi(t)|, \ \forall t \in \left[t_0 - \bar{h}, t_0\right]$ □

To simplify further referencing, we present a homogeneous counterpart of (3.10) as follows,
$$\dot{y} = p(t)y + c(t)L\big(t, y(t), y(t - h_1(t)),..., y(t - h_m(t))\big)$$
$$y(t,\phi) = |\varphi(t)|, \ \forall t \in \left[t_0 - \bar{h}, t_0\right] \quad (3.12)$$

and assume that (3.12) admits a unique solution under the conditions of Theorem 1 as well.

Next, let us assume that Definitions 1 and 2 are concerned with the solutions to equations (3.12) and (3.10) and that $r_i$ are the superior values of $\delta_i$, i.e., $r_i \equiv \sup \delta_i$, $i = 1,...,7$, for which the $i$-th condition in Definition 1 holds for (3.12) or the conditions 6 or 7 in Definition 2 hold for (3.10). Consequently, $r_1 = r_1(t_0)$, $r_3 = r_3(t_0)$, $r_6 = r_6(t_0, F_o), \ \forall F_0 \leq F_*(t_0)$ and $r_7 = r_7(F_o), \ \forall F_0 \leq \hat{F}$. Then let $B_{r_i} := \{\varphi \in C\big(\left[t_0 - \bar{h}, t_0\right], \mathbb{R}^n\big): \ \|\varphi\| \leq r_i\}$, $i = 1,...,7$ be the balls with radiuses $r_i$ that are centered at the origin.

This embraces the following statements.

**Theorem 2**. Assume that the conditions of Theorem 1 are met and the trivial solution to equation (3.12) is either stable, uniformly stable, asymptotically stable, uniformly asymptotically stable, or exponentially stable. Then the trivial solution to equation (2.15) is, in turn, stable, uniformly stable, asymptotically stable, uniformly asymptotically stable, or exponentially stable, respectively.

Furthermore, under the conditions of 1- 5 of Definition 1, $B_{r_i} \subseteq B_{R_i}$, $i = 1,...,5$, respectively.

**Proof**. Both above statements directly follow from the application of inequality (3.11) to the solutions to equations (2.15) and (3.12) □

**Theorem 3**. Assume that the conditions of Theorem 1 are met and that the solutions to the equation (3.10) are bounded for the set values of $t_0$ and $F_0 \leq F_*(t_0)$ if $\|\varphi\| \leq r_6$ or uniformly bounded for $F_0 \leq \hat{F}$ if $\|\varphi\| \leq r_7$. Then the solutions to equation (2.14) with the matched history functions are also bounded or uniformly bounded for the same values of $t_0$ and $F_0$, respectively. Furthermore, under the prior conditions $B_{r_i} \subseteq B_{R_i}$, $i = 6,7$.

**Proof**. The proof of this statement directly follows from the application of inequality (3.11) to the solutions to equation (2.14) and its scalar counterpart (3.10) □

The key task of estimating the values of $r_i$, $i = 1,...,7$ can be abridged through utility of the following statement.

**Lemma 5**. Assume that $y_i(t,\phi_i), \ \forall t \geq t_0$, are the solutions to (3.10) or (3.12) and $\phi_1(t) \geq \phi_2(t), \forall t \in \left[t_0 - \bar{h}, t_0\right]$. Then $y_1(t,\phi_1) \geq y_2(t,\phi_2), \ \forall t \geq t_0$.

**Proof**. The proof of this statement follows from the assumption of uniqueness of solutions to (3.10)/(3.12) which prohibits the intersection of the solution curves to a scalar delay equation in $t \times y$ - plane □

Additionally, the proof of Lemma 5 can also be fostered by the application of Lemma 4 under the assumption that $L(t, \zeta_1,...,\zeta_{m+1})$ is a nondecreasing function in all its variables starting from the third, which is often met in applications, see Appendix.

Thus, a solution to equations (3.10) or (3.12) with a constant history function, i.e., $y = y(t,\phi), \ \forall t \geq t_0$, $\phi(t) = q \in \mathbb{R}, \ \forall t \in \left[t_0 - \bar{h}, t_0\right]$, monotonically increases in $q$, $\forall t \geq t_0$ which further streamlines the simulations of the values of $r_i$.

Next, we apply our approach to ensure the robust stability of the trivial solution to the vector delay equation (2.15) through the conditions applied to its scalar counterpart. Consequently, we write a scalar counterpart for

(2.17) as follows

$$\dot{z} = p(t)z + c(t)L(t, z(t), z(t-h_1(t)),...,z(t-h_m(t))) + c(t)L_R(t, z(t), z(t-h_1^*(t)),...,z(t-h_m^*(t))), \forall t \geq t_0$$
$$z(t,\phi) = |\varphi(t)|, \forall t \in [t_0 - \bar{h}, t_0]$$
(3.13)

where $z \in \mathbb{R}_{\geq 0}$ and $L_R(t, \xi_1,...,\xi_{m+1}) \in \mathbb{R}_{\geq 0}$, $\forall \xi_i \in \mathbb{R}_{\geq 0}$, $\forall t \geq t_0$ is a continuous function in all its variables. Then the application of Definition 3 to (3.12) and (3.13) leads to the following statement.

**Theorem 4**. Assume that conditions of Theorem 1 and conditions (I)-(IV) hold for equation (3.13) and that the trivial solution to (3.12) is uniformly asymptotically stable. Then the trivial solution to (2.15) is robustly stable.

**Proof**. In fact, under the above conditions, the trivial solution of (3.12) is robustly stable which implies that $z(t,|\varphi|) < \varepsilon$, $\forall t \geq t_0$, where $z(t,|\varphi|)$ is a solution to (3.13), if $|L_R(t, \xi_1,...,\xi_{m+1})|_\infty < \Delta_1(\varepsilon)$, $\forall t \geq t_0$, $\forall \xi_i < \varepsilon$, $i = 1,...,m+1$, $\|\varphi\| < \Delta_2(\varepsilon)$, and $|h_i(t) - h_i^*(t)| < \Delta_3(\varepsilon)$, $\forall t \geq t_0$, see Definition 3. Then (3.11) yields that $|x(t,\varphi)| \leq z(t,|\varphi|) < \varepsilon$, $\forall t \geq t_0$, where $x(t,\varphi)$ is a solution to (2.17), since (3.13) is the scalar counterpart to (2.17) □

Lastly, we show that the FTS's conditions for equations (2.14) and (2.15) can be applied to their scalar counterparts as well. This leads to the following statement.

**Theorem 5**. Assume that:
  a) equation (3.10) or (3.12) is FTS with respect $\alpha, \beta, T \in \mathbb{R}_+$, $\alpha < \beta$. Then equations (2.14) or (2.15) is FTS with respect $\alpha, \beta, T$ as well.
  b) equation (3.10) or (3.12) is FTCS with respect to $\alpha, \beta, \gamma, T \in \mathbb{R}_+$, $\alpha < \gamma$. Then equations (2.14) or (2.15) is FTCS with respect $\alpha, \beta, \gamma, T$ as well.

**Proof**. In fact, both above statements directly follow from (3.11) □

Thus, the developed approach embraces the assessment of the boundedness/stability characteristics of nonautonomous VNDS through streamlined simulations of their scalar counterparts which also estimate the radiuses of the balls encompassing history functions stemming bounded or stable solutions to the original system.

## 4. A Closed-Form Robust Stability Criterion

An analytical criterion for boundedness or stability, initially developed for a scalar delay equation (see, e.g., [6] and additional references), can often be extended to the corresponding vector system using the methodology outlined above. Moreover, simplified versions of equations (3.10) and (3.12) can be obtained by bounding their right-hand sides with time-invariant or linear functions. These approaches will be demonstrated in this and subsequent sections.

The following statement provides a criterion for the robust stability of the trivial solution to a nonlinear and non-autonomous system of ODEs, accounting for nonlinear perturbations that involve variable delays and coefficients.

**Corollary 1**. Assume that the conditions of Theorem 1 hold and $f \in C([t_0, \infty) \times \mathbb{R}^n; \mathbb{R}^n)$, $\sup_{\forall t \geq t_0} p(t) = \hat{p} < 0$, $\sup_{\forall t \geq t_0} c(t) = \hat{c} < \infty$, $L \in C([t_0, \infty) \times \mathbb{R}_{\geq 0}; \mathbb{R}_{\geq 0})$, $\sup_{\forall t \geq t_0} L(t, y) = \hat{L}(y) < \infty$. Next, we also presume that, $\hat{p} y + \hat{c}\hat{L}(y) < 0$, $\forall y \in (0, y_+)$, $y_+ > 0$, and stated above conditions (I)-(IV) hold for (3.13). Then the trivial solution to (2.15) is robustly stable.

**Proof**. In fact, due to the set conditions, the trivial solution to the scalar equation (3.12) is uniformly asymptotically stable. Then the application of Theorem 4 ensures this statement □

## 5. Linearized Delay Auxiliary Equations

This section linearizes the scalar equations (3.10) or (3.12) near the origin of these systems by applying the Lipschitz-like continuity condition. This approach simplifies the criteria for the boundedness and stability of the scalar equations, and consequently, of their vector counterparts.

Let us assume that,

$$L(t, \zeta_1,...,\zeta_{m+1}) \leq \sum_{i=1}^{m+1} \mu_i(t, \tilde{\zeta})\zeta_i, \forall |\zeta| \leq \tilde{\zeta} > 0, \forall t \geq t_0$$
(5.1)

where $\zeta_i \in \mathbb{R}_{\geq 0}$, $\zeta = [\zeta_1,...,\zeta_{m+1}]^T$ and $\mu_i \in C([t_0,\infty) \times \mathbb{R}_{\geq 0}; \mathbb{R}_{\geq 0})$. Then the application of (5.1) to (3.10) yields the following scalar linear equation,

$$\dot{u} = P(t)u + c(t)\left(\sum_{i=2}^{m+1} \mu_i(t,\tilde{\zeta})u(t-h_i(t)) + F_0|e(t)|\right) \quad (5.2)$$
$$u(t,\phi) = |\varphi(t)|, \forall t \in [t_0 - \overline{h}, t_0]$$

where $P(t) = p(t) + c(t)\mu_1(t,\tilde{\zeta})$. To shorten further citations, we also write a homogeneous counterpart to (5.2),

$$\dot{u} = P(t)u + c(t)\sum_{i=2}^{m+1} \mu_i(t,\tilde{\zeta})u(t-h_i(t)) \quad (5.3)$$
$$u(t,\phi) = |\varphi(t)|, \forall t \in [t_0 - \overline{h}, t_0]$$

Next, we set in (5.1) that,
$$\zeta_1 = u(t,\phi), \zeta_i = u(t-h_i(t),\phi) \quad (5.4)$$

where $u(t,\phi)$ is either a solution to equation (5.2) or equation (5.3). In the latter case, we can conclude that inequality (5.1) holds if (5.3) is stable and $\|\varphi\|$ is sufficiently small. Similarly, in the former case, we demonstrate that inequality (5.1) holds if the preceding conditions are met and $F_0$ is also sufficiently small.

These observations lead to the following statements.

**Lemma 6**. Assume that $\mu_i \in C([t_0,\infty) \times \mathbb{R}_{\geq 0}; \mathbb{R}_{\geq 0})$, $i=1,...,m+1$, $p \in C([t_0,\infty); \mathbb{R})$, $c \in C([t_0,\infty); \mathbb{R}_+)$, $\varphi \in C([t_0,\infty); \mathbb{R}^n)$, functions $h_i(t)$ comply with (2.4), equation (5.3) admits a unique solution $\forall \tilde{\zeta} < \zeta_* > 0$, $\forall t \geq t_0$ and $\forall \|\varphi\| \leq \hat{\varphi} > 0$. Also, we assume that the trivial solution to (5.3) is stable for the set value of $t_0$ and $\forall \tilde{\zeta} \leq \tilde{\zeta}_{max} < \zeta_*$, where $\tilde{\zeta}_{max}$ is the maximal value of $\tilde{\zeta}$ for which (5.3) is stable. Then, inequality (5.1) holds if $\|\varphi\|$ is a sufficiently small, where $\zeta_i(t,\phi)$ are assigned by (5.4) and $u(t,\phi)$ is a solution to (5.3).

**Proof**. In fact, Definition 1a implies that under these conditions, $u(t,\phi) < \varepsilon$, $\forall t \geq t_0$, where $\varepsilon$ is an arbitrary small value and $\phi(t) = |\varphi(t)|$, $\forall t \in [t_0 - \overline{h}, t_0]$ □

In turn, this leads to the following stability criteria.

**Theorem 6**. Assume that the conditions of Theorem 1 and Lemma 6 hold. Then,
$$|x(t,\varphi)| \leq y(t,\phi) \leq u(t,\phi), \forall t \geq t_0 \quad (5.5)$$
Furthermore, the trivial solutions to both equations (3.12) and (2.15) are stable.

**Proof**. Really, under these conditions (5.1) holds due to Lemma 6. In turn, the application of Lemma 4 yields that $y(t,\phi) \leq u(t,\phi)$, $\forall t \geq t_0$, which subsequently implies stability of the trivial solutions to (3.12) and (2.15) □

**Corollary 2**. Assume that the conditions of Theorem 6 hold and, also, the trivial solution to (5.3) is either uniformly stable, asymptotically stable, uniformly asymptotically stable, or exponentially stable. Then, the trivial solutions to equations (3.12) and (2.15) will exhibit the same type of stability: either uniformly stable, asymptotically stable, uniformly asymptotically stable, or exponentially stable, respectively.

**Proof**. In fact, inequality (5.1) and (5.5) hold under all the above conditions which ensures this statement □

To further simplify equation (5.3), we assume that $\sup_{\forall t \geq t_0} p(t) = \hat{p} < 0$, $\sup_{\forall t \geq t_0} c(t) = \hat{c} < \infty$, $\sup_{\forall t \geq t_0} \mu_i(t) = \hat{\mu}_i < \infty$ and review (5.3) as follows,

$$\dot{U} = PU + \hat{c}\left(\sum_{i=2}^{m+1} \mu_i(\tilde{\chi})U(t-h_i(t))\right), U(t,\phi) = |\varphi(t)|, \forall t \in [t_0 - \overline{h}, t_0] \quad (5.6)$$

where $P = \hat{p} + \hat{c}\hat{\mu}_1$ and $U \in C([t_0,\infty); \mathbb{R}_{\geq 0})$. This leads to the following statement.

**Corollary 3**. Assume that (5.6) admits a unique solution $\forall \tilde{\zeta} < \zeta_\times \leq \zeta_*$, $\forall t \geq t_0$ and $\forall \|\varphi\| \leq \varphi_\times \leq \hat{\varphi}$. Also assume that (5.6) is stable for the set value of $t_0$ and $\forall \tilde{\zeta} \leq \tilde{\zeta}_U \leq \tilde{\zeta}_{max}$, and $\|\varphi(t)\|$ is a sufficiently small number. Then,

$$|x(t,|\varphi|)| \le y(t,|\varphi|) \le u(t,|\varphi|) \le U(t,|\varphi|), \forall t \ge t_0$$

and the trivial solutions to equations (5.3), (3.12) and (2.15) are stable.

**Proof**. The proof of this statement imitates the proofs of Lemma 6 and Theorem 6 □

Thus, the stability of linear scalar delay equation (5.3) and its autonomous counterpart (5.6) imply stability of the trivial solution to a vector system (2.15). Some stability criteria for two former equations, for instance, can be found in [7], see also additional references therein. Note also that the stability criteria for (5.6) are readily known if $h_i \equiv const$, see, for instance, [13],[16], [23].

In turn, we extend this approach to establish sufficient boundedness criteria for equation (2.14) by analyzing the solutions of the linear scalar equation (5.2) which can be represented as follows, see, e.g., [23][23], p.139,

$$u(t,\phi) = u_h(t,\phi) + F_0 u_{nh}(t,0) \tag{5.7}$$

where $u_h(t,\phi)$ is a solution to (5.3) and $u_{nh}(t,0)$ is a solution to (5.2) with $\varphi(t) \equiv 0$ and $F_0 = 1$. In turn,

$$u_{nh}(t,0) = \int_{t_0}^{t} C(t,s)c(s)|e(s)|ds$$, where $C(t,s)$ is a Cauchy function for (5.2) which is defined as a solution to (5.3)

with the following history function, $C(t,s) = \begin{cases} 0, & t_0 - s_1 \le t < s \\ 1, & t = s \end{cases}$. Thus, $C(t,s) \ge 0$ which yields that

$u_{nh}(t,0) \ge 0, \forall t \ge t_0$. This leads to the following statement.

**Lemma 7**. Assume that the conditions of Lemma 6 are met, and that equation (5.2) has a unique solution under these conditions. Also, we assume that $|u_{nh}(t,0)| \le \hat{u}_{nh} < \infty, \forall t \ge t_0$, and $F_0$ and $\|\varphi\|$ are sufficiently small. Then, inequality (5.1) holds, where $\zeta_i(t,\phi)$ are assigned by (5.4) and $u(t,\phi)$ is a solution to (5.2).

**Proof**. In fact, under the above conditions, both additions on the right side of (5.7) can be made arbitrary small $\forall t \ge t_0$ by the appropriate choice of both $F_0$ and $\|\varphi(t)\|$ □

Next, we introduce a criterion for the boundedness of solutions to vector equation (2.14), derived from a condition imposed on the scalar equation (3.2).

**Theorem 7**. Assume that the conditions of Theorem 1 and Lemma 7 are met, and both $F_0$ and $\|\varphi\|$ are sufficiently small. Then,

$$|x(t,\varphi)| \le y(t,\phi) \le u(t,\phi), \forall t \ge t_0 \tag{5.8}$$

where $x(t,\varphi)$, $y(t,\phi)$ and $u(t,\phi)$ are solutions to equations (2.14), (3.10) and (5.2), respectively, and $\phi(t) = |\varphi(t)|, \forall t \in \left[t_0 - \bar{h}, t_0\right]$.

**Proof**. In fact, (5.8) immediately follows from the application of Lemma 7 and Lemma 4 □

Thus, (5.7) and (5.8) acknowledge that,

$$|x(t,\varphi)| \le u_h(t,|\varphi(t)|) + F_0 u_{nh}(t,0) \tag{5.9}$$

where, under conditions of Theorem 7, $u_{nh}(t,0) < \infty$ and $u_h(t,\phi)$ is stable, asymptotically stable, or exponentially stable solution to (5.3) depending on the applicable criteria. Thus, (5.9) establishes the input-to-state stability of the vector delay equation (2.14) based on the conditions applied to its scalar and linear counterpart. For additional references on input-to-state stability for delay systems, see [3].

The last two sections introduce further simplified linear auxiliary equations that provide more conservative but easier-to-apply criteria for assessing the boundedness and stability of systems (2.14) and (2.15). However, a more accurate evaluation of the dynamics of the original vector systems can be achieved through straightforward simulations of their scalar counterparts (3.10) or (3.12). These simulations not only validate our methodology but also gauge its accuracy.

## 6. Simulations

Let us assume that (2.14) takes the following form,

$$\dot{x} = (A_0(t) + A_1(t))x + \rho A_1(t)x(t-h) + f(t, x(t-h)) + F(t)$$
$$x(t) \in \mathbb{R}^2, x(t) = x_0, \forall t \in \left[-\bar{h}, 0\right] \tag{6.1}$$

where $f(t,x(t-h)) = b[0, x_2^3(t-h)]^T$, $F(t) = [0, F_2]^T$, $b, h, \bar{h}, \rho \in \mathbb{R}_+$, $x_0 \in \mathbb{R}^2$ is a constant vector,

$$A_0(t) = \begin{pmatrix} \lambda(t) & 0 \\ 0 & \lambda(t) \end{pmatrix}, \quad A_1(t) = \begin{pmatrix} 0 & 1 \\ -\omega(t) & -\alpha_1 \end{pmatrix}.$$

Next, we calculate functions $p(t)$ and $c(t)$ for the subsystem $\dot{x} = A_0(t)x$, see Section 3. Clearly, in this case, $w(t) = diag(\eta, \eta)$, where $\eta = \exp\left(\int_{t_0}^{t} \lambda(s)ds\right)$, which implies that $p(t) = \lambda(t)$ and $c(t) \equiv 1$. Thus, the scalar auxiliary aquation (3.10) for (5.1) takes the following form,

$$\dot{y} = (\lambda(t) + |A_1(t)|)y + \rho|A_1(t)|y(t-h) + |b|y^3(t-h) + |F(t)|, \quad y(t) \in \mathbb{R}_{\geq 0}$$
$$y(t) = |x_0|, \quad \forall t \in [-h_1, 0] \tag{6.2}$$

In turn, an autonomous scalar counterpart to (6.2) can be written as follows,

$$\dot{\hat{y}} = (\hat{\lambda} + \hat{A}_1)\hat{y} + \rho\hat{A}_1\hat{y}(t-h) + |b|\hat{y}^3(t-h) + F_0, \quad \hat{y}(t) \in \mathbb{R}_{\geq 0}$$
$$\hat{y}(t) = |x_0|, \quad \forall t \in [-h, 0] \tag{6.3}$$

where $\sup_{\forall t \geq t_0} \lambda(t) = \hat{\lambda} < 0$, $\sup_{\forall t \geq t_0} |A_1(t)| = \|A_1\| = \hat{A}_1$. Thus, due to Lemma 4, we infer that $y(t) \leq \hat{y}(t)$, $\forall t \geq 0$.

Let us recall that the stability of the linear equation $\dot{\hat{y}} = (\hat{\lambda} + \hat{A}_1)\hat{y} + \rho\hat{A}_1\hat{y}(t-h)$ is well studied in the space of its parameters and the asymptotic stability of this equation implies the asymptotic stability of the trivial solution to the homogeneous counterpart to (6.3) if $|b|$ is sufficiently small [23]. Furthermore, the right sides of both equations (6.2) and (6.3) depend only on a scalar function $|A_1(t)|$ and scalar parameter $\hat{A}_1$, rather than matrix $A_1(t)$.

To describe the results of simulations of equations (6.1)-(6.3), we set below that $\lambda = \lambda_0 + \lambda_+(t)$, where $\lambda_0 = -3$ and either a) $\lambda_+(t) = q\sin(dt)$, $q = 0.1$, $d = 5$ or b) $\lambda_+(t) = q\exp(-dt)$, $q = 1$, $d = 1$. Note that in the former case $\hat{\lambda} = -2.9$ whereas in the later one $\hat{\lambda} = -2$. Next, we also assume that $\omega(t) = \omega_0 + \omega_1(t)$, $\omega_1(t) = a_1\sin(r_1t) + a_2\sin(r_2t)$, where $\omega_0 = 1$, $a_1 = a_2 = 0.1$, $r_1 = 1$, $r_2 = 3.14$. Additionally, we set in our simulations that $e(t) = [0, \sin(r_3t)]^T$, $r_3 = 10$; the values of other parameters are indicated below in the captions of the figures. All simulations presented in this paper utilize the MATLAB code DDE23.

Figure 1 contrasts in solid, dashed, and dashed dotted lines the norms of solutions to (6.1) and the corresponding solutions to (6.2) and (6.3), simulating with matching history functions. All plots demonstrate that $|x(t,x_0)| \leq y(t,|x_0|) \leq \hat{y}(t,|x_0|)$ on the simulated time interval which validates our theoretical inferences.

Furthermore, the solutions to scalar nonautonomous equation (6.2) provide fairly accurate approximations to the norms of the matched solutions to (6.1) across a wide range of parameters that were considered in our simulations. As expected, simulations of the autonomous equation (6.3) deliver more crude estimates. Additionally, the accuracies of these estimates are more sensitive than the former ones to enlargements of the delay parameter $h$, the scale of the delay nonlinear term $b$ and the amplitude of input $F_0$. Nevertheless, the solutions to the scalar nonautonomous equation (6.2) still adequately estimate the norms of the matched solutions to (6.1) in the wide range of parameters considered in our simulations.

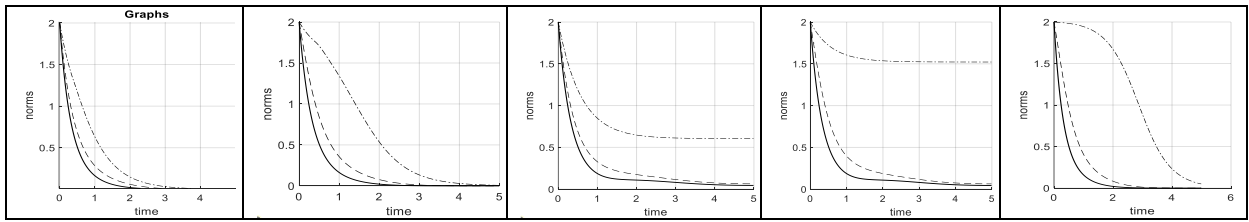

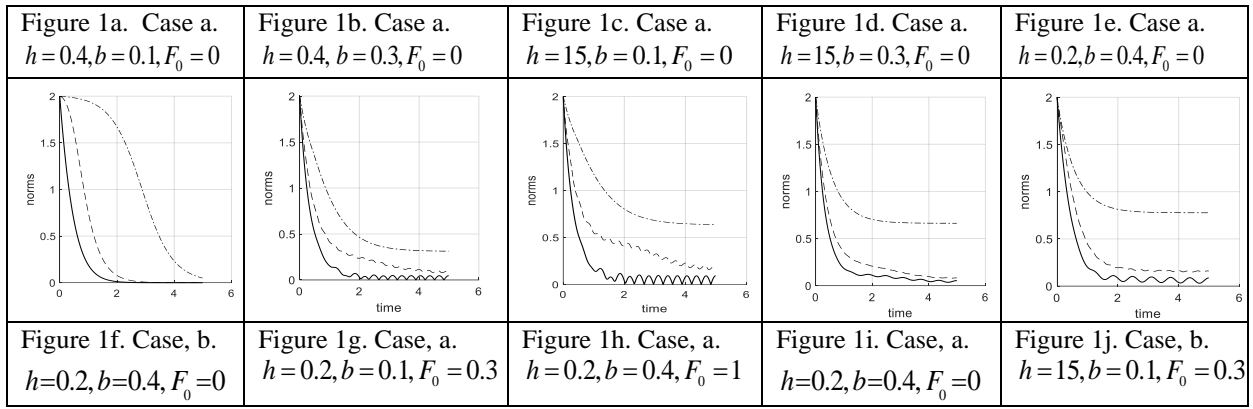

| Figure 1a. Case a. $h=0.4, b=0.1, F_0=0$ | Figure 1b. Case a. $h=0.4, b=0.3, F_0=0$ | Figure 1c. Case a. $h=15, b=0.1, F_0=0$ | Figure 1d. Case a. $h=15, b=0.3, F_0=0$ | Figure 1e. Case a. $h=0.2, b=0.4, F_0=0$ |
|---|---|---|---|---|
| Figure 1f. Case, b. $h=0.2, b=0.4, F_0=0$ | Figure 1g. Case, a. $h=0.2, b=0.1, F_0=0.3$ | Figure 1h. Case, a. $h=0.2, b=0.4, F_0=1$ | Figure 1i. Case, a. $h=0.2, b=0.4, F_0=0$ | Figure 1j. Case, b. $h=15, b=0.1, F_0=0.3$ |

Figure 1. Comparison of time evolution of the norms of solutions to equation (6.1) and the corresponding solutions to scalar equations (6.2) and (6.3), simulated with matching history functions.

To simulate the boundaries of the trapping/stability regions, which include the initial vector entries leading to bounded or stable solutions to equation (6.1) or its homogeneous counterpart, we represent $x_0$ in the polar coordinates as, $x_{0,1} = r_0 \cos\varphi_0$, $x_{0,2} = r_0 \sin\varphi_0$. The angle coordinate was then discretized with the step size equals $\pi/100$. For each angle, the radial coordinate was adjusted sequentially to approximate the values located at the boundary of the region. The simulations employed a binary search, starting from a small radial value, which was progressively increased until a rapid growth in $|x(t, x_0)|$ was numerically observed.

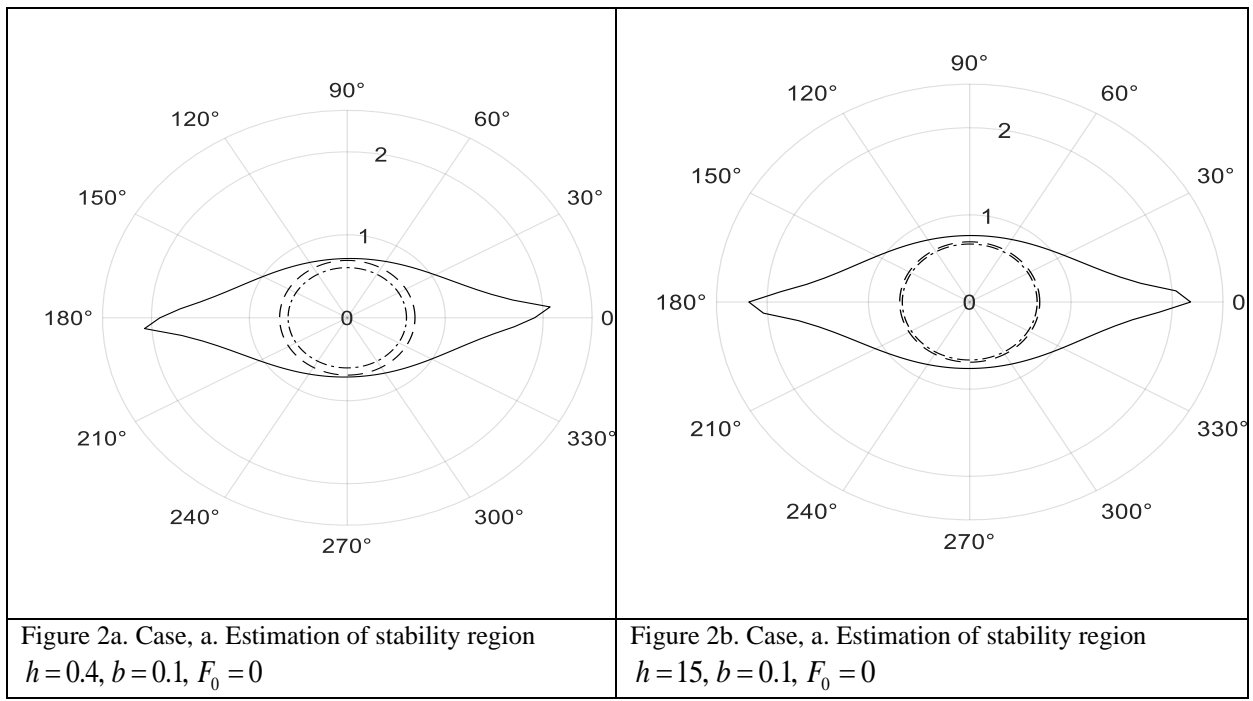

| Figure 2a. Case, a. Estimation of stability region $h=0.4, b=0.1, F_0=0$ | Figure 2b. Case, a. Estimation of stability region $h=15, b=0.1, F_0=0$ |

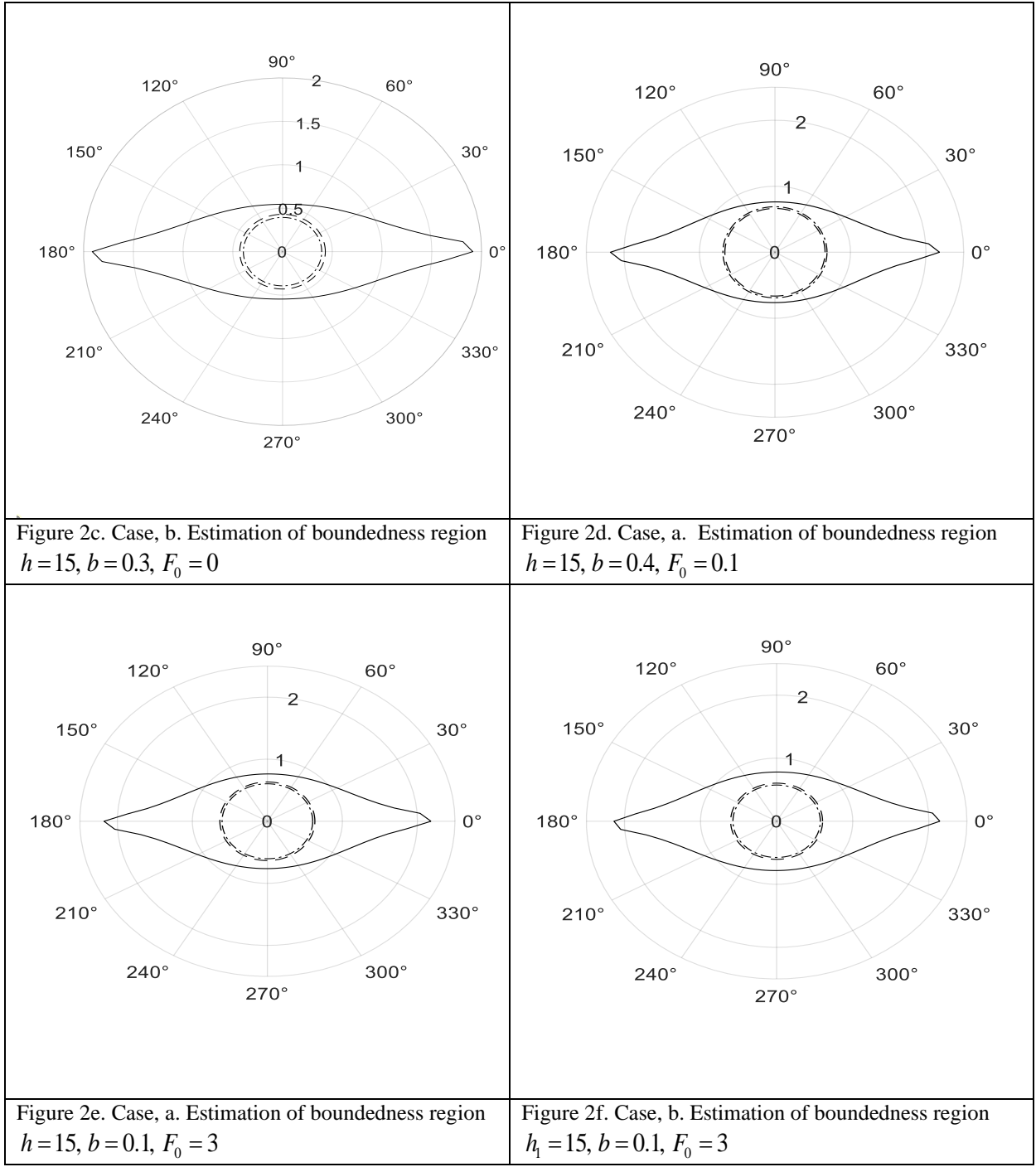

| Figure 2c. Case, b. Estimation of boundedness region $h=15, b=0.3, F_0=0$ | Figure 2d. Case, a. Estimation of boundedness region $h=15, b=0.4, F_0=0.1$ |
|---|---|
| Figure 2e. Case, a. Estimation of boundedness region $h=15, b=0.1, F_0=3$ | Figure 2f. Case, b. Estimation of boundedness region $h_1=15, b=0.1, F_0=3$ |

Figure 2. Estimation of the radii of disks embedded within trapping/stability regions of equation (6.1) or its homogeneous counterpart.

A similar approach was applied to scalar equations (6.2) and (6.3) to approximate the radii of disks embedded within the trapping/stability regions of equation (6.1) or its homogeneous counterpart. The time required to estimate these radiuses is practically insensitive to enlargement of the number of equations coupled in the original vector system and is further reduced since both $y(t,|x_0|)$ and $\hat{y}(t,|x_0|)$ monotonically increase in $|x_0|$ for $\forall t \geq t_0$. In contrast, the time required to estimate the boundary of the trapping/stability regions for the original equation (2.14) is scaled by $O(m^n)$, where $m$ is the number of points that were taken to discretize a phase-space variable and $n$ is the number of these variables.

Figure 2 shows in solid, dashed and dashed dotted lines the boundaries of the trapping ($F_0 \neq 0$) or stability ($F_0 = 0$) regions that were obtained in simulations of equations (6.1), (6.2), and (6.3), respectively, for the different sets of parameters noted in the figure's captions. The curves in Figure 2 are plotted in polar coordinates, where the natural logarithm of radius vector is plotted for each curve instead of actual radius.

The simulations of equations (6.2) and (6.3) are used merely to estimate the radii of disks within the trapping/stability regions of (6.1) or its homogeneous counterpart. Interestingly, these estimates are consistently close to each other across all plots in Figure 2 and align well with the smallest dimensions of the trapping/stability regions. Furthermore, these estimates of the radii remain largely unaffected by changes in the delay parameter, the amplitude of the input, or the scale of the nonlinear component. Finally, the arrangement of the curves in all plots in Figure 2 supports our theoretical inferences.

These straightforward estimations of the radii of disks within the boundedness/stability regions of vector nonlinear systems with delay and variable coefficients are particularly valuable in practical applications. This is especially true in cases where the system's history function is unknown, but the range of variation of its superior norm can still be reasonably gauged.

## 7. Conclusion

This paper presents a novel approach that develops scalar counterparts for a broad class of vector nonlinear systems with variable coefficients and multiple varying delays. The solutions to these scalar equations serve as an upper bound for the norms of the solutions to the original equations, provided that the history functions of the corresponding scalar and vector equations are properly matched. This enables the evaluation of the boundedness and stability characteristics of vector delay nonlinear systems by analyzing their scalar counterparts, which can be done through seamless simulations or simplified analytical methods. As a result, we derived new boundedness and stability criteria and estimated the radii of the balls that encompass history functions leading to bounded or stable solutions for the original vector system.

Additionally, we demonstrate that applying a Lipschitz-like condition to our nonlinear scalar equations, reveals their linear counterparts. The solutions of these scalar linear delay equations provide upper bounds for the solutions of their nonlinear scalar counterparts, further simplifying analytical inferences and leading to additional boundedness and stability criteria.

Furthermore, the developed approach aggregates parameters, history functions, and models of uncertainties, simplifying the assessment of the original systems' robustness properties.

Our results are validated through representative simulations, showing that the solutions to the scalar auxiliary equations fairly approximate the norms of solutions to the original vector systems if the latter are stemmed from the central regions of the boundedness/stability domains. Furthermore, we demonstrated that our method estimates the radii of the boundedness/stability balls of the test system with relatively small error.

Future work will focus on extending the applicability of this methodology and developing efficient recursive approximations for the boundaries of the trapping/stability regions and bilateral solution bounds for the norms of solutions to vector nonlinear and non-autonomous systems with delay.

**Appendix**. Let us show how to derive inequality (3.3) using a simple example that can be naturally extended to more complex cases. Assume that $x = [x_1 \ x_2]^T$ and a vector-function $f$ is defined, e.g., as follows, $f = [a_1(t)x_1^3 x_2^2(t-h_1) \quad a_2(t)x_2^3(t-h_2)]^T$. Then $|f|_2 \leq |f|_1 \leq |a_1||x_1^3||x_2^2(t-h_1)| + |a_2||x_2^3(t-h_2)|$ $\leq |a_1||x(t)|^3|x(t-h_1)|^2 + |a_2||x(t-h_2)|^3$, where we use that $|x_i^n| \leq |x|^n$, $|x_i^n(t-h)| \leq |x(t-h)|^n$, $i=1,2$, $n \in \mathbb{N}$. Clearly, this inference can be extended to the power series and some rational functions under the pertained conditions.

**Acknowledgement**. The code for simulations discussed in Section 6 was developed by Steve Koblik
Data sharing is not applicable to this article as no datasets were generated or analyzed during the current study.
The author declares that he has no conflict of interest.